\documentclass{amsart}

\usepackage{epsfig}
\usepackage{amsthm}
\usepackage{amssymb}
\usepackage{amsmath}
\usepackage{amscd}

\newtheorem {Theorem}   {Theorem} 

\numberwithin{Theorem}{section}

\theoremstyle{definition}
\newtheorem{Definition}[Theorem]{Definition}
\theoremstyle{remark}
\newtheorem{Remark}[Theorem]{Remark}

\newtheorem{Claim}[Theorem]{Claim}

\expandafter\chardef\csname pre amssym.def at\endcsname=\the\catcode`\@ 
\catcode`\@=11 
\def\undefine#1{\let#1\undefined} 
\def\newsymbol#1#2#3#4#5{\let\next@\relax 
 \ifnum#2=\@ne\let\next@\msafam@\else 
 \ifnum#2=\tw@\let\next@\msbfam@\fi\fi 
 \mathchardef#1="#3\next@#4#5}
\def\mathhexbox@#1#2#3{\relax 
 \ifmmode\mathpalette{}{\m@th\mathchar"#1#2#3}%
 \else\leavevmode\hbox{$\m@th\mathchar"#1#2#3$}\fi} 
\def\hexnumber@#1{\ifcase#1 0\or 1\or 2\or 3\or 4\or 5\or 6\or 7\or 8\or 
 9\or A\or B\or C\or D\or E\or F\fi} 
 
\font\teneufm=eufm10 
\font\seveneufm=eufm7
\font\fiveeufm=eufm5 
\newfam\eufmfam
\textfont\eufmfam=\teneufm 
\scriptfont\eufmfam=\seveneufm 
\scriptscriptfont\eufmfam=\fiveeufm 
 
\catcode`\@=\csname pre amssym.def at\endcsname 

\def	\eps	{\epsilon}

\def	\C      {{\mathbb C}}
\def	\R	{{\mathbb R}}

\def	\Q	{{\mathbb Q}}

\def	\T	{{\mathbb T}}

\def    \ra     {{\rightarrow}}
\def    \haf    {{\frac{1}{2}}}

\begin{document}






\title[new counterexamples]{New smooth counterexamples to the Hamiltonian Seifert Conjecture}

\author[Ely Kerman]{Ely Kerman}
\address{The Fields Institute, 222 College Street, Toronto, Ontario M5T 3J1, Canada}
\email{ ekerman@fields.utoronto.ca}

\date{\today}

\thanks{This work was supported by an NSERC postdoctoral fellowship }


\bigskip
  
\begin{abstract}
We construct a new aperiodic symplectic plug and hence new smooth counterexamples to the Hamiltonian Seifert conjecture in $\R^{2n}$ for $n\geq 3$. In other words, we develop an alternative procedure, to those of V. L. Ginzburg \cite{gi1,gi2} and M. Herman \cite{her}, for constructing smooth Hamiltonian flows, on the standard symplectic $\R^{2n}$ for $n\geq3$, which have compact regular level sets that contain no periodic orbits. The plug described here is a modification of those built by Ginzburg. In particular, we utilize a different ``trap'' which makes the necessary embeddings of this plug much easier to construct.
\end{abstract}

\maketitle

\section{Introduction} 

The problem of establishing the existence of periodic orbits on fixed level sets of Hamiltonian flows has been studied extensively and there are now many theorems which achieve this for a variety of level sets and symplectic manifolds, see e.g. \cite{ch,fhv,ho,hv,ke,vi}. Helping to determine the limits of these existence theorems is a small set of examples of Hamiltonian flows with compact regular level sets containing no periodic orbits. The most important of these examples are defined on the standard symplectic $\R^{2n}$ for $n\geq 3$. These were constructed, independently and in different ways, by V. L. Ginzburg \cite{gi1,gi2} and M. Herman \cite{her}. More precisely, smooth Hamiltonian flows with aperiodic level sets were built for $n>3$ in \cite{gi1} and \cite{her}; a  $C^{3-\eps}$-example for $n=3$ also appeared in \cite{her}; and, finally, a $C^{\infty}$-example for $n=3$ was built in \cite{gi2}. For a thorough discussion of the relevance of these examples to the various existence theorems, as well as a list of all the (previously) known examples of Hamiltonian flows without periodic orbits, the reader is referred to \cite{gi3} and \cite{gi4}. 

In this paper, we describe new examples of smooth aperiodic Hamiltonian flows in $\R^{2n}$ for $n\geq3$. These are built using Ginzburg's symplectic version of a method due to Wilson \cite{wi}, which uses aperiodic plugs to create aperiodic flows \cite{gi1}. In particular, we build a new symplectic plug (see Definition \ref{def:1}) which is similar to the ones constructed in \cite{gi1} and \cite{gi2} but which utilizes a different trap to destroy periodic orbits. The advantage of the plug described here is that the flow on the trap is quite simple and the necessary embeddings of the plug into Euclidean space are much easier to find. 

In the remainder of this section, we briefly recall Wilson's procedure for producing aperiodic flows and discuss how it has been used to find counterexamples to the Seifert conjecture. In the second section, we describe Ginzburg's extension of this method to Hamiltonian flows. We then construct our new symplectic plug in the third section of the paper.

\subsection{Counterexamples to the Seifert and Hamiltonian Seifert conjectures }  One of the most well known conjectures concerning the existence of periodic orbits is the Seifert conjecture. This is the assertion that every nonvanishing vector field on $S^3$ must have a periodic orbit. Since it was posed (asked) in \cite{se}, three counterexamples to this conjecture have been constructed, first in the $C^1$-category by P. Schweitzer \cite{sch}, then in the $C^{3-\epsilon}$-category by J. Harrison \cite{ha}, and ultimately in the $C^{\omega}$-category by K. Kuperberg \cite{kku1}. For a review of these constructions see \cite{kku2}.

One can extend this conjecture to Hamiltonian flows in a number of ways \cite{gi3}. For the purposes of this paper, the Hamiltonian Seifert conjecture will be the assertion that every compact regular level set of a smooth function on a symplectic manifold must contain a periodic orbit of the corresponding Hamiltonian flow. Just as its namesake, this conjecture is known to be false. In particular, the constructions of Ginzburg, Herman and those defined here, can be used to produce counterexamples for any symplectic manifold of dimension at least six. However, the conjecture remains open for $\R^{4}$. It also continues to merit study because the rare examples of aperiodic Hamiltonian flows may help to illuminate the poorly understood boundary between existence and nonexistence results, (see $\S 3.4$ of \cite{gi4}). 

All of the counterexamples mentioned above are built using a procedure which was introduced by Wilson in \cite{wi}. One begins with a manifold $W$ and a vector field ${X}$, both in the appropriate class, such that the periodic orbits of ${X}$ are well understood. Then local alterations are made to ${X}$ in order to open these closed trajectories. More precisely, suppose that ${X}$ has a periodic orbit which passes through the point $w \in W$. We would like to alter $X$ near $w$ so that 
\begin{enumerate}
\item for the resulting vector field, the trajectory through $w$ is no longer closed
\item no new periodic orbits are created
\item the new vector field is still in the desired class.
\end{enumerate}
To perform such alterations Wilson introduced the  notion of an aperiodic plug. This consists of a compact manifold $P$ (with boundary), a carefully constructed vector field $V$ on $P$ and an embedding $j\colon P \to \R^{n}$ (where $n=\dim W$) such that $dj(V)$ is constant near $\partial P$. With this, the dynamical system defined by $V$ on $P$ can be inserted, via $j$, into a flow box around $w$. This insertion is made in such a way that the trajectories near $w$ become trajectories of $V$ whose flow has been tailored to accomplish the goals above. For example, with the exception of K. Kuperberg's plug, the first goal is achieved by forcing the trajectory through $w$ to asymptotically approach an aperiodic orbit on an invariant submanifold $M \subset P$. This trajectory is hence trapped in the image of $P$ and because of its asymptotic behavior is no longer closed. The submanifold $M$ together with its flow is called the trap (or core) of the plug. 

In \cite{wi}, Wilson uses the irrational flow on a two-dimensional torus as a trap and constructs a smooth plug which can be embedded in a flow box of dimension at least four. Schweitzer then constructs a $C^1$-plug in \cite{sch} which can be embedded in a flow box of dimension three and thus obtains the first counterexample to the original Seifert conjecture. The trap which allows Schweitzer to do this is a punctured two-dimensional torus with a flow that is trivial near the boundary and contains the (aperiodic) closure of an (aperiodic) orbit of the Denjoy flow in the interior. Unlike the complete torus, the punctured torus can be embedded into $\R^3$ in such a way that the standard normal vectors are  everywhere parallel to a constant vector. This fact is what enables Schweitzer to suitably embed his plug into $\R^3$ and his crucial observation can be summarized as follows:
\begin{quote} The flow on the trap only needs to be aperiodic on an invariant closed subset. \end{quote} The strength of this simple statement lies in the fact that we know of very few examples of totally aperiodic flows on closed manifolds of small dimension. In other words, it greatly increases the set of flows which may be used as the trap of a plug.  

\begin{Remark}
A remarkable feature of K. Kuperberg's plug is that, instead of using a known aperiodic flow as a trap, she creates an entirely new aperiodic dynamical system in dimension three (see \cite{kku1,gh}).    
\end{Remark}

In \cite{gi1} and \cite{gi2}, Ginzburg uses a trap consisting of the horocycle flow on $ST^*\Sigma$, the unit cotangent bundle of a surface $\Sigma$ with constant negative curvature. This flow is not only aperiodic on all of $ST^*\Sigma$, but is also minimal \cite{hed}. In order to use this trap to obtain the counterexample for $n=3$, it was necessary for Ginzburg to improve, in \cite{gi2}, the dimensional constraints for some symplectic embedding results of M. Gromov \cite{gr1,gr2}.  Utilizing Schweitzer's observation above, our choice of a trap is a relatively simple flow on $\T^3$ which is only aperiodic on an invariant subset. For this trap, the construction of the necessary embeddings is much simpler. Unfortunately, this choice does {\bf not} allow us to reduce to the case $n=2$. Most likely, the construction of a counterexample for this final case (if it is even possible!) will require, as in \cite{kku1}, the creation of an entirely new type of Hamiltonian dynamics. 

\subsection*{Acknowledgments.}I would like to express my gratitude to Viktor Ginzburg for his many encouraging and helpful comments. In particular, the relative simplicity of the construction presented here owes much to one of his suggestions.

\section{ Symplectic plugs and aperiodic Hamiltonian flows }

In this section we recall how the special nature of a Hamiltonian flow allows one to reformulate the Hamiltonian Seifert conjecture and to use a different type of plug to produce counterexamples.

We begin by considering the existence of periodic orbits for a nonvanishing vector field $X$ on an odd-dimensional manifold $W$. In the presence of an additional geometric structure on $W$ this question can be rephrased. Recall that a two-form is said to be maximally nondegenerate if its kernel is everywhere one-dimensional. Let us assume that we can find such a form $\sigma$ on $W$ so that $X(w)$ spans $\ker \sigma(w)$ for all $w\in W$. Then, instead of looking for periodic orbits of $X$, we may equivalently look for smooth maps $\gamma \colon S^1 \to W$ for which $\dot{\gamma}(t)$ spans $\ker \sigma (\gamma (t))$ for all $t \in [0,1].$ Such maps are called closed characteristics of $\sigma$ on $W$ and we will identify those having the same image. If $\gamma$ is instead defined on some interval, i.e. the image isn't a loop, it will just be called a characteristic of $\sigma$.   

Now, let $H$ be a smooth function on a symplectic manifold $(Q,\Omega$) and let $W$ be a compact regular level set of $H$. Denote by $X_H$ the Hamiltonian vector field defined by the equation $$i_{X_H}\Omega = dH.$$ This vector field is tangent to $W$ and nonvanishing there. It is also easy to verify that $\sigma = \Omega |_{W}$ is maximally nondegenerate and that $X_H$ spans $\ker \sigma$. From the discussion above, we see that the Hamiltonian Seifert conjecture can be reformulated as a statement which asserts the existence of closed characteristics for certain maximally nondegenerate two-forms. In particular, the conjecture concerns those forms which are obtained by  restricting a symplectic form to a hypersurface in a symplectic manifold. With this formulation of the conjecture in mind, Ginzburg introduced the following notion of a symplectic plug in \cite{gi1} as a tool to construct counterexamples.  

\begin{Definition}
\label{def:1} A {\em symplectic plug} for dimension $2n$  $(n>2)$ is a product manifold $P= [-1,1] \times N^{2n-2}$, equipped with an exact maximally nondegenerate two-form $\omega$ and embeddings $j \colon P \to \R^{2n-1}$ and $J \colon P \to \R^{2n}$ such that the following conditions hold:

\begin{description}
\item[(P.1)] $\ker{\omega}$ is vertical (parallel to $[-1,1]$) near $\partial P.$
\item[(P.2)] there exists a point $(-1, \tilde{p}) \in \{-1\} \times N$ such that any characteristic of $\omega$ through this point does not pass through $\{1\} \times N.$
\item[(P.3)] a characteristic through $(-1,p) \in \{-1\} \times N$ which meets $\{1\} \times N$ must do so at $(1,p)$.
\item[(P.4)] $\omega$ has no closed characteristics in $P$.
\item[(P.5)] $j^* \Omega_{2n-2}= \omega$ near $\partial P,$ where $\Omega_{2n-2}$ is the pullback of the canonical symplectic form on $\R^{2n-2}$ via the projection $\R^{2n-1} \to \R^{2n-2}$.
\item[(P.6)] for a given $\delta >0$ we have $J \colon P \to [-\delta, \delta] \times\R^{2n-1}$ and $J^* \Omega_{2n}= \omega$. Moreover, if we let $\tilde{j}$ be the map $j$ with its range identified with  $\{0\} \times \R^{2n-1}\subset \R^{2n}$, then  $J=\tilde{j}$ near $\partial P$ and $J$ is isotopic to $\tilde{j}$ relative to $\partial P$.
\end{description}
\end{Definition}

\begin{Remark} 
\label{rmk:f}
It is perhaps unnecessary to include the embedding $j$ in this definition. However, as described below, this map allows a symplectic plug to be used for the slightly different purpose of destroying closed characteristics for any maximally nondegenerate two-form. It also helps to clarify the construction of the embedding $J$. Indeed, here and in \cite{gi1} the map $J$ is of the form $J= f \circ j$ where $f \colon \R^{2n-1} \to \R^{2n}.$ While $j$ is again used to destroy a closed characteristic, the map $f$ now ensures that this is done without disturbing the symplectic structure. When the domain of $f$ is identified with $\{0\} \times \R^{2n-1} \subset \R^{2n}$ we will denote it by $\tilde{f}$ so that $J= \tilde{f}   \circ \tilde{j}$.\end{Remark} 

Symplectic plugs were constructed for $n \geq 3$ in \cite{gi1} and \cite{gi2}. In \cite{gi1}, Ginzburg also describes how a symplectic plug for dimension $2n$ can be used to produce counterexamples to the Hamiltonian Seifert conjecture in $\R^{2n}$. This procedure is described below.

As a first step, we will show that a symplectic plug can be used to destroy closed characteristics of a maximally nondegenerate two-form $\sigma$ on an odd-dimensional manifold  $W^{2n-1}$. In fact, assuming that $\sigma$ has only a finite number of closed characteristics we will construct from it a new maximally nondegenerate two-form with none. Let $w \in W$ be a point in the image of one of the closed characteristics of $\sigma$ and choose a small open neighborhood of $w$ which does not meet the others. This neighborhood can be chosen to be diffeomorphic to an open ball $B^{2n-1} \subset \R^{2n-1}$ in such a way that $\sigma$ gets identified with $\Omega_{2n-2}$. This is a version of a flow box around $w$ since the characteristics of $\Omega_{2n-2}$ in $\R^{2n-1}$ are all straight. We can then use $j$ to insert $(P, \omega)$ into this neighborhood so that $(-1, \tilde{p})$ gets mapped to $w$. Let $\sigma^{\prime}$ be the form which equals $\sigma$ outside $j(P)$ and $\omega$ inside. By (P.1) and (P.5),  we see that $\sigma^{\prime}$ is well-defined, smooth and maximally nondegenerate. Properties (P.2), (P.3) and (P.4) then imply that the characteristics of $\sigma^{\prime}$ through $w$ are no longer closed and $\sigma^{\prime}$ has no new closed characteristics. We have thus succeeded in finding a maximally nondegenerate two-form on $W$ with exactly one less closed characteristic than $\sigma$. Applying this procedure to the other closed characteristics we have shown that the existence of symplectic plugs for $n \geq 3$ implies the following theorem.
\begin{Theorem}[\cite{gi1, gi2}]
\label{thm:maxnondeg}
Assume that $\dim W \geq 5$ and $\sigma$ has a finite number of closed characteristics. Then there exists a closed maximally nondegenerate two-form $\sigma^{\prime}$ which is homotopic to $\sigma$ and has no closed characteristics.
\end{Theorem}
\begin{Remark}
\label{rmk:hom}
Two closed maximally nondegenerate two-forms are said to be homotopic if they can be joined by a one parameter family of such forms all lying in the same cohomology class. The existence of the homotopy in Theorem \ref{thm:maxnondeg} follows from one between $\omega$ and a two-form $\omega_1$ which has a vertical kernel and satisfies $j^* \Omega_{2n-2} = \omega_1$. This homotopy is defined in Remark \ref{rmk:homo} and is not a necessary feature for a symplectic plug.  
\end{Remark}
We can now use the embedding $J$ to extend the previous theorem to one which applies to a hypersurface $W$ in a symplectic manifold $(Q^{2n},\Omega).$ Assume that $\Omega|_W$ has a finite number of closed characteristics and let $w \in W$ be a point on one of them. There exists a neighborhood $U$ of $w$ in $Q$ which doesn't meet the other closed characteristics on $W$ and is also a symplectic flow box, i.e. it is diffeomorphic to an open ball $B^{2n} \subset \R^{2n} = \R \times \R^{2n-1}$ in such a way that the intersection of $B^{2n}$ with $\{0\} \times \R^{2n-1}$ corresponds to $U \cap W$, $\Omega$ gets identified with $\Omega_{2n}$ and $\Omega|_{U \cap W}$ gets identified with $\Omega_{2n-2}$ on ${\{0\} \times  \R^{2n-1}}$. We then use $J$ (corresponding to a sufficiently small $\delta$) to insert $(P,\omega)$ into $U$ so that $J(\partial P) \subset U \cap W $ and $(-1, \tilde{p})$ gets mapped to $w$. Let $W^{\prime}$ be the hypersurface which is equal to $W$ outside of $U$ and agrees with the image of J inside of $U$. We claim that $\Omega|_{W^{\prime}}$ has exactly one less closed characteristic than $\Omega|_W$. To see this, we view the insertion of $(P, \omega)$, by $J= \tilde{f} \circ \tilde{j}$, as a two step process (see Remark \ref{rmk:f}). The primary insertion uses $\tilde{j}$ to map $(P, \omega)$ into $U \cap W$. As above, this yields a new form $\Omega^{\prime}$ on $W$ with exactly one less closed characteristic than $\Omega|_{W}$. The map $\tilde{f}$ then deforms $W$ into $W^{\prime}$ so that  $\Omega^{\prime}={\tilde{f}}^* \Omega$. In particular, $\tilde{f}$ takes the  characteristics of $\Omega^{\prime}$ through $w \in W$ to the characteristics of $\Omega|_{W^{\prime}}$ through $w \in W^{\prime}$, which can therefore never be closed. Repeating this process for the other closed characteristics we get;
\begin{Theorem}[\cite{gi1, gi2}]
\label{thm:hyp}  
Assume that $\dim Q \geq 6$ and that $\Omega|_W$ has a finite number of closed characteristics. Then there exists a smooth hypersurface $W^{\prime} \subset Q$, $C^0$-close and isotopic to $W$, such that $\Omega|_{W^{\prime}}$ has no closed characteristics. 
\end{Theorem} Letting $W$ be an irrational ellipsoid in $(\R^{2n}, \Omega_{2n})$ with $n \geq 3$, the latter result yields concrete examples of smooth hypersurfaces (regular level sets) without closed characteristics (periodic orbits), see Corollary 3.3 and Corollary 3.4 in \cite{gi3}.
 
\section{A new symplectic plug}
We will now describe the construction of a new symplectic plug for $n=3$. This procedure can easily be generalized to work for $n>3$. 

\subsection{The trap}
As a trap, we will use the constant-speed geodesic flow of the flat metric on the two-dimensional torus. In terms of a Hamiltonian flow, this is defined on $(T^*\T^2, d\lambda)$, the cotangent bundle of the two-dimensional torus equipped with its canonical symplectic structure. With respect to the (global) coordinates $\{\theta_1,\theta_2, p_1, p_2 \} $ on  $T^*\T^2= \T^2 \times \R^2$, let $H$ be the standard kinetic energy function defined by 
\begin{align*}
H \colon  T^*\T^2 &\to \R\\
          (\theta_1,\theta_2,p_1,p_2) &\mapsto \haf(p_1^2+p_2^2).
\end{align*}
Set $M$ equal to $ \{H=R^2\}$ and denote by $S$ the dense subset of $M$ consisting of points with $(p_1,p_2)$ satisfying $\frac{p_2}{p_1} \notin \Q$. It is easy to check that the Hamiltonian vector field of $H$ on $M$ has only base components and that the flow on the torus $\T^2 \times \{(p_1,p_2)\}$ is the constant linear flow with slope $\frac{p_2}{p_1}$. Hence, each torus corresponding to a rationally independent pair $(p_1,p_2)$ is an invariant closed subset of $M$ on which the flow is aperiodic. Alternatively, this situation can be described as follows. The manifold $M$ is diffeomorphic to the three dimensional torus and is equipped with an exact, maximally nondegenerate two-form $d\eta = d\lambda|_M$. Any characteristic which starts in $S$ stays in $S$ and all such characteristics are not closed. The manifold $M$ together with the form $d\eta$ will be the trap of our symplectic plug.   

Switching to polar coordinates $(p_1,p_2) \to (r,\theta_3)$ we see that $\bar{\theta}=(\theta_1,\theta_2,\theta_3)$ is a set of coordinates on $M$ and $$\eta = R\cos(\theta_3) d\theta_1 + R\sin(\theta_3) d\theta_2.$$ We will call the angle $\theta_3$ irrational if it corresponds to a rationally independent pair $(p_1,p_2)$.
    
\subsection{The plug}

Let $N= [-\eps,\eps]\times M$ so that $P=[-1,1]\times [-\eps,\eps]\times M$. We label the coordinates for the two intervals by $t$ and $x$, respectively. Consider the two-form on $P$ defined by $$\omega = d(A(\theta_3,x,t) \eta + B(x,t)dt),$$ for smooth functions $A$ and $B$ of the indicated coordinates.

\begin{Claim}Let the functions $A$ and $B$ have the following properties:
\begin{description}
\item [(A.1)] $A>0$.
\item [(A.2)] $A^{\prime}_x\geq0$, with equality only at $(\tilde{\theta}_3,0,\pm \haf)$ where $\tilde{\theta}_3$ is a fixed irrational angle. These are also the only two critical points of $A$.
\item [(A.3)] $A=1+x$ near $\partial P$.
\item [(A.4)] $A$ is even in $t$.
\item [(B.1)] $B=0$ near $\partial P$.
\item [(B.2)] $B=x$ near $(0,\pm \haf)$.
\item [(B.3)] $B$ is odd  in $t$.
\end{description}
Then, $\omega$ is maximally nondegenerate on $P$ and (P.1), (P.2), (P.3) and (P.4) hold.
\end{Claim}

\begin{Remark} Functions satisfying the conditions in the claim are easily constructed.
\end{Remark}

\begin{proof}

Treating $\omega$ as a family of skew--symmetric matrices parameterized by $P$, it is easy to check that the coefficient of the linear term in the characteristic polynomial of $\omega(\theta_1,\theta_2,\theta_3,x,t)$ is $$R^4 A^2({A^{\prime}_x}^2+{A^{\prime}_t}^2) + R^2{B^{\prime}_x}^2( A^2 + {A^{\prime}_{\theta_3}}^2). $$ Clearly, $\omega$ is maximally nondegenerate if and only if this coefficient is never zero. By properties (A.1) and (A.2), the first term only vanishes at $(\tilde{\theta}_3,0,\pm \haf)$. At these points the second term becomes $$ R^2 {B^{\prime}_x}^2 A^2$$ which does not equal zero by property (B.2). Thus, $\omega$ is maximally nondegenerate.

Near the boundary of $P$ we have $$\omega = dx \wedge \eta + (1+x)d\eta.$$ Since $\eta$ is nondegenerate on $\ker d\eta$, we see that $\frac{\partial}{\partial t}$ spans $\ker{\omega}$ near $\partial P$. Hence, $\omega$ has property (P.1).

In order to check (P.2), (P.3) and (P.4), we consider a nonvanishing vector field $V$ on $P$ such that $V$ spans $\ker{\omega}$. We choose $V$ so that it equals $\frac{\partial}{\partial t}$ near $\partial P$  and we recall that any characteristic of $\omega$ is equivalent (after some reparameterization) to a trajectory of $V$.

It is straightforward to check that $$\omega \wedge \omega \wedge dt = R^2A A^{\prime}_x \, \mu,$$ where $\mu$ is the standard volume form with respect to our coordinates on $P$. Since $i_V \omega=0$, it follows that the $t$-component of $V$, $V_t$, satisfies $$V_t \,\omega \wedge \omega =  R^2A A^{\prime}_x \, i_V \mu.$$ This equation implies that $V_t$ vanishes only when $R^2A A^{\prime}_x$ does. Hence, by property (A.2), $V_t$ equals zero along the two-dimensional tori $\T^2_{\pm}$ defined by $(\theta_3,x,t)= (\tilde{\theta}_3,0,\pm \haf)$, respectively. By continuity and our choice of $V$ near $\partial P$, we also see that $V_t$ is strictly positive away from these tori. On $\T^2_{\pm}$, $$ \omega = A(\tilde{\theta}_3,0,\pm\haf)d\eta + dx \wedge dt.$$ Hence, $V$ is tangent to both tori and lies in the kernel of $d \eta$ on each. This means that the flows of $V$ on $\T^2_{\pm}$ are both conjugate to the irrational flow for the angle $\tilde{\theta}_3$. Let $L^-(\T^2_-)$ be the negative limit set of $\T^2_-$ under the flow of $V$. The previous observations then imply that all the points in  $$\left(\{-1\}\times N\right) \cap L^-(\T^2_-)$$ asymptotically approach the aperiodic flow on $\T^2_-.$ Consequently, the trajectories of $V$ through all of these points are trapped in $P$ and property (P.2) holds. 

At this point, it is also easy to see that none of the trajectories of $V$ are closed and hence (P.4) is satisfied. Since $V_t$ is nonnegative, any periodic orbit of $V$ would have to lie on either $\T^2_{+}$ or $\T^2_{-}$. However, this is impossible since the flow of $V$ on both these sets is aperiodic. 

Now, let $\psi \colon P \to P$ be the map which sends $t \to -t$ and acts as the identity on the $N$-component. By properties (A.4) and (B.3), it is possible to choose $V$ such that $$d\psi (V) = -V.$$ Consider a trajectory of $V$ (a characteristic of $\omega$) which starts at some point $(-1,p) \in \{-1\} \times N$ and exits $P$ through $\{1\} \times N$. The $t$-component $V_t$ can not vanish along this trajectory since it would then be trapped in one of the internal tori $\T^2_{\pm}$ and could never exit $P$. Consequently, the (anti-)symmetry of $V$ described above implies that as the trajectory progresses from $t=0$ to $t=1$ it retraces, in reverse, the progress it made along the $N$-component from $t=-1$ to $t=0$. Hence, the trajectory must exit at $(1,p)$ and property (P.3) has been verified.
\end{proof}

In order to complete the construction of our plug we must find embeddings $j\colon P \to \R^{5}$ and $J\colon P \to \R^{6}$ which satisfy (P.5) and (P.6). First, consider the restriction of $\omega$ to $$\{-1\}\times N \approx [-\eps,\eps]\times M.$$ We call this restriction $\rho$ and note that $$\rho = dx \wedge \eta + (1+x)d \eta. $$ This form restricts to the hypersurface $\{0\}\times M$ as $d\eta.$ Notice that $[-\eps,\eps]\times M$ can also be considered as a (closed) neighborhood of $M$ in  $T^*\T^2$ and recall that the canonical symplectic form $d\lambda$ also restricts to $M$ as $d\eta$. Using Weinstein's Extension theorem \cite{we}, it is easy to show that two symplectic forms which agree on a compact oriented hypersurface are equivalent in a neighborhood of that hypersurface, see \cite[Exercise 3.35, p.102]{ms}. Hence, for sufficiently small $\eps$ there exists a symplectomorphism $$\phi_1 \colon ([-\eps,\eps]\times M, \rho) \to (\mathcal{N}(M), d\lambda),$$ where $\mathcal{N}(M)$ is an open neighborhood of $M$ in $T^*\T^2$. Next, for sufficiently small $R$ (and $\eps$), we can assume that $\mathcal{N}(M)$ lies in an arbitrarily small neighborhood of the zero section in $T^*\T^2$. Invoking the Lagrangian Neighborhood theorem, there exists a symplectomorphism $$ \phi_2 \colon (\mathcal{N}(M), d\lambda) \to (\R^4, \Omega_4)$$ onto a neighborhood surrounding a Clifford torus in $\R^4$. Composing $ \phi_1 $ and $ \phi_2$ we get a symplectomorphism $$\phi \colon ([-\eps,\eps]\times M, \rho) \to (\R^4, \Omega_4).$$ We now define the embedding $j$ by 
\begin{align*}
j \colon  P=[-1,1] \times N &\to \R \times\R^{4}=\R^{5} \\
          (t,p) &\mapsto (t,\phi(p)).
\end{align*}
Since $\omega = \rho$ near all of $\partial P$, $j$ satisfies (P.5).

\begin{Remark}\label{rmk:homo} To form the homotopy discussed in Remark \ref{rmk:hom}, we replace $A$ and $B$ in the definition of $\omega$ by the family of functions
\begin{eqnarray*}
A_{\tau}&=& (1-\tau)A + \tau(1+x)\\ 
B_{\tau}&=& (1-\tau)B.
 \end{eqnarray*}
The resulting two-forms, $\omega_{\tau}$, are all maximally nondegenerate and by the construction of $j$ we also have $j^*\Omega_{4}=\omega_1$.    
\end{Remark}
As noted in Remark \ref{rmk:f}, the embedding $J$ is of the form $J=f\circ j$ for some embedding $f \colon j(P)\subset \R^5 \to \R^6$. For simplicity we will identify $j(P)$ with $P$. We then choose $f$ so that it maps $P$ into $\R \times P \subset \R \times \R^5$. On $\R \times P$ we have coordinates $(y,t,x,\bar{\theta})$. Recalling that $j^*(\Omega_4)=\rho$ we see that with respect to these coordinates
\begin{eqnarray*}
{\Omega_6|}_{\R \times P} &=& dy \wedge dt + \rho\\ 
         &=& d( ydt + (1+x)\eta).
\end{eqnarray*}
Let $f \colon P \to \R \times P $ be defined by $$(t,x,\bar{\theta}) \mapsto (B, t, A-1,\bar{\theta}).$$ Property (A.2) ensures that $f$ is indeed an embedding. In addition, we may choose the function $B$ to be arbitrarily small so that given any $\delta >0$ the image of $f$ will lie in $(-\delta, \delta)\times \R^5$.  Note also that $A$ and $B$ are equal to $1+x$ and $0$ on $\partial P$, respectively, and are also isotopic to these functions. This implies that $f$ equals $id$ near the $\partial P$ and is isotopic to $id$ relative  to $\partial P$. Hence, we only need to prove that $f^* \Omega_6 =\omega$. In fact, we can verify the stronger fact that $$f^*(\underbrace{ydt + (1+x)\eta}_{\alpha})= A\eta +Bdt.$$ To see this let $v \in T_pP.$ Then 
\begin{eqnarray*}
f^* \alpha (v) &=& \alpha(f)[df(v)]\\
&=& (A\eta + Bdt)[df(v)]\\
&=& (A\eta + Bdt)[v],
\end{eqnarray*}
where the last equality holds because $f$ acts like the identity on $t$ and $\bar{\theta}$.
With this, the construction of our symplectic plug is complete.


\begin{thebibliography}{ABCD}

\bibitem[Ch]{ch}
Chen, W., Pseudo-holomorphic curves and the Weinstein conjecture, {\em Comm. Anal. Geom.}, {\bf 8}, 115-131 (2000).

\bibitem[FHV]{fhv}
Floer, A., Hofer H., Viterbo, C., The Weinstein conjecture in $P \times \C^l$, {\em Math. Z.}, {\bf 203}, 469--482 (1990).

\bibitem[Gh]{gh}
Ghys, E., Construction de champs de vecteurs sans orbite p\'{e}riodique (d'apr\`{e}s Krystyna Kuperberg), 
{\em S\'{e}m. Bourbaki}, 1993--1994, Ast\'{e}risque, no 227, Soc. Math. France, Paris, 283--307 (1994).

\bibitem[Gi1]{gi1}
Ginzburg, V.L., An embedding $S^{2n-1}\ra \R^{2n},\,2n-1 \geq 7$ whose Hamiltonian flow has no periodic trajectories, {\em IMRN}, {\bf 2}, 83-97 (1995).

\bibitem[Gi2]{gi2}
Ginzburg, V.L., A smooth counterexample to the Hamiltonian Seifert conjecture in $\R^6$, {\em IMRN}, {\bf 13}, 641-650 (1997).

\bibitem[Gi3]{gi3}
Ginzburg, V.L., Hamiltonian dynamical systems without periodic orbits, {\em Amer. Math Soc. Transl.}, (2) {\bf 196}, 35-48 (1999).

\bibitem[Gi4]{gi4}
Ginzburg, V.L., The Hamiltonian Seifert conjecture: examples and open problems, to appear in {\em Proceedings of the Third European Congress of Mathematics}, Preprint 2000; math.DG/0004020.


\bibitem[Gr1]{gr1}
Gromov, M., A topological technique for the construction of solutions of differential equations and inequalities, {\em Proc. Internat. Congr. Math. (Nice, 1970)}, Vol.2, Gauthier-Villars, Paris, 221-225 (1971).

\bibitem[Gr2]{gr2}
Gromov, M., {\em Partial differential relations}, Springer-Verlag, New York, 1986.


\bibitem[Ha]{ha}
Harrison, J., A $C^2$ counterexample to the Seifert conjecture, {\em Topology}, {\bf 27}, 249-278 (1988).

\bibitem[Hed]{hed}
Hedlund, G.A., Fuschian groups and transitive horocycles, {\em Duke Math. J.}, {\bf 2}, 530-542 (1936).

\bibitem[Her]{her}
Herman, M.-R., Examples of compact hypersurfaces in $\R^{2p}$, $2p \geq6$, with no periodic orbits, in {\em Hamiltonian systems with three or more degrees of freedom}, C. Simo (Editor), NATO Adv. Sci. Inst. Ser. C, Math. Phys. Sci., vol. 533, Kluwer Acad. Publ., Dordrecht, (1999).

\bibitem[Ho]{ho}
Hofer, H., Pseudoholomorphic curves in symplectizations with applications to the Weinstein conjecture in dimension three, \emph{Invent. Math.}, {\bf 114}, 515-563 (1993). 

\bibitem[HV]{hv}
Hofer, H., Viterbo, C., The Weinstein conjecture for cotangent bundles
and related results, \emph{Ann. Scuola Norm. Sup. Pisa Cl. Sci.},
(4) {\bf 15} (1988) no. 3, 411--445 (1989).


\bibitem[HZ]{hz}
Hofer, H., Zehnder, E., Periodic solution on hypersurfaces and a
result by C. Viterbo, {\em Invent. Math.}, {\bf 90}, 1--9 (1987).

\bibitem[Ke]{ke}
Kerman, E., Periodic orbits of Hamiltonian flows near symplectic critical submanifolds, {\em IMRN}, {\bf 17}, 953-969 (1999).

\bibitem[KuG]{gku}
Kuperberg, G., A volume-preserving counterexample to the Seifert conjecture, \emph{Comment. Math Helv.}, {\bf 71}, 239-268 (1996).

\bibitem[KuK1]{kku1}
Kuperberg, K., A smooth counterexample to the Seifert conjecture in dimension three, \emph{Ann. of Math.}, (2) {\bf140}, 723-732 (1994).  

\bibitem[KuK2]{kku2}
Kuperberg, K., Counterexamples to the Seifert conjecture, \emph{Proc. Internat. Congr. Math. (Berlin 1998)}, Vol. II,  831-840.  



\bibitem[MS]{ms}
McDuff, D., and Salamon, D., \emph{Introduction to symplectic topology},
Oxford University Press, New York, (1995).


\bibitem[Sc]{sch}
Schweitzer, P.A., Counterexamples to the Seifert conjecture and opening closed leaves of foliations, \emph{Ann. of Math.}, (2) {\bf100}, 229-234 (1970).

\bibitem[Se]{se}
Seifert, H., Closed integral curves in 3-space and isotopic two-dimensional deformations, {\em Proc. Amer. Math. Soc.}, {\bf1}, 287-302 (1950).

\bibitem[Vi]{vi}
Viterbo, C., A proof of the Weinstein conjecture in $\R^{2n}$, {\em Ann. Inst. H. Poincar\'{e}, Anal. Non Lin\'{e}are}, {\bf 4}, 337-356 (1987).

\bibitem[We]{we}
Weinstein, A., {\em Lectures on symplectic manifolds}, volume 29 of CBMS, {\em Reg. Conf. Ser. in Math.} AMS, 1977.


\bibitem[Wi]{wi}
Wilson, F., On the minimal sets of non-singular vector fields, \emph{Ann. of Math.}, (2) {\bf84}, 529-536 (1966).
\end{thebibliography}
\end{document}